# Higher-order social-ecological network as a simplicial complex

Sudeepto Bhattacharya

Department of Mathematics
and
Centre for Environmental Sciences and Engineering

School of Natural Sciences, Shiv Nadar University Delhi-NCR, Gautam Buddha Nagar, Greater Noida 201 314 Uttar Pradesh, India

e-mail: sudeepto.bhattacharya@snu.edu.in

**Abstract**

A social-ecological network is a formal representation of a corresponding social-ecological system, and encodes a relation within a given system as an interaction. Conventionally, such networks have been defined as encoding and representing pairwise interactions among the fundamental units of the system. This work proposes a combinatorial definition of social-ecological network by means of its structure as a simplicial complex. The proposed definition is a comprehensive one that takes into account the heterogeneity of interactions within a given SES, and the higher-order social-ecological network modelled using this definition is able to represent the modelled SES by capturing all orders of interactions within the system. Such a social-ecological network consequently, is better equipped to capture and represent the structural details of the real-world SES, and is thus capable of facilitating a deeper insight into the complex behaviour of the represented SES emergent through the higher-order interactions within the system, as compared to the conventional graph-theoretic network that exclusively models pairwise interactions.

**Introduction**

A social-ecological system (SES) is a complex adaptive system [1, 2, 3]. For every such system, there exist corresponding social-ecological networks (SEN) that encode (that is, convert the information into a mathematical form) the interactions in the system, giving rise to a formal representations of the complex SES. SENs as network theoretic representations have offered a very powerful and an equally successful framework for obtaining insight into the inherent complexity of the modelled SES, along with an understanding of its various behaviour and properties, and over the past two decades have evolved as an epitome of formal modelling of SESs [4 – 17].

Most of the SENs proposed and studied in the scholarship, however, are premised on the framework of simple graphs and therefore are based on the assumption that the complex interactions and structural behaviour and related properties of the SES can be modelled exclusively by considering the binary or pairwise interactions among the units of the system. Such a graph theory-based SEN formalism for a SES completely models the represented system provided the system admits only pairwise interactions among its fundamental units. However, for SESs that may admit not only binary but also non-pairwise, higher-order interactions, such SEN representations face a twin mathematical challenge:(i) due to the definition of graph as a combinatorial object, the SEN model is necessarily limited to only capture the binary relations existing in the system, and hence are successful in capturing the system interactions only partially, and therefore, (ii) there exists a possibility of masking those complex behaviours of the modelled system that are expressed due to non-pairwise, higher-order interactions in the system. These two limitations often cause a risk of producing an oversimplified network model which fails to capture the complexity arising due to higher-order interactions in the real system, and hence fails to be a faithful (in the sense of being a close approximation) representation of the modelled SES. The principal reason for this failure lies in the definition of SEN, which is defined using the graph theoretic framework, allowing only for binary, or pairwise interactions to be encoded. It is therefore of an overriding need that mathematical frameworks beyond the graph theoretic one be adopted to define SEN so that it is structurally able to capture and represent all orders of interactions in the corresponding SES.

The objective of this work is to construct a definition for SEN that is a structurally generalised version of the conventional graph theoretic definition, using the methods of combinatorial algebraic topology. The proposed definition leads to a framework that addresses both the problems mentioned in the earlier paragraph satisfactorily. While the proposed definition of SEN serves as a generic definition for all networks – random, regular or complex - it takes into consideration not only the existence of pairwise interactions among the system units as conventional networks do, but also the usually more abundant higher-order interactions within the SES, which remain excluded in the graph theory-based SENs, as discussed above.

The next section provides the necessary mathematical preliminaries required to achieve the objective of this work. Similar definitions as presented in this section are found elsewhere, and have been repeated here to make the work self-contained. We begin with the definition of structure as a mathematical object, essentially following the detailed discussions presented in texts such as [18, 19, 20]. Next, we employ this definition to formalise the concept of an SES through different orders of interactions among its constituent fundamental units with the objective to arrive at a definition of an SEN representing it.

**Preliminaries**

*Structure*

*Definition*: A structure is a tuple $\theta = (X, R, F, C)$ where:

    (i)    $X$ is a non-empty set, the universe or the domain or the carrier set of the structure;

    (ii)    $R = \{R_1,...,R_m\}$ is the set of relations on $X$;

    (iii)    $F = \{F_1,...,F_n\}$ is the set of mappings on $X$;

    (iv)    $C = \{c_i \mid i \in I\}$, $c_i \in X \;\; \forall i \in I$ is the set of constants, where $I$ is an index set.

The *signature* of the structure is a sequence $(F_1, ..., F_n; R_1, ..., R_m; \kappa)$, where $F_i : X^{f_i} \to X$ is an $n$-ary map for each $i$; $R_j \subseteq X^{r_j}$, is an $m$-ary relation for each $j$; $\kappa = |\{c_i \mid i \in I\}|$, cardinality of $I$, on $X$ where $n, m \in \mathbb{N}$.

Often, it is the universe $X$ that forms the focus of study, as all the mappings, relations and constants are defined on $X$. As such, the structure may also be written as $\theta_X = (X, F, R, C)$, or more commonly as $X = (X, F, R, C)$, thus identifying the structure with the universe. Further, a structure in which the set $F = \emptyset$ is called a relational structure, and would form the central component of what follows in this work. Sometimes, however, different symbols may be used when there is a need to distinguish between and emphasise on the name of the object defined by the structure and the universe of the object. We shall adopt the above notational conventions interchangeably and as per need in what follows, as we apply the definition to formally define a system.

*System*

*Definition*: A *system* $\Sigma$ is an object given by the structure $\theta_\Sigma = (V, E, C)$, where $V$, the *vertex set* is a non-empty set, the universe of the system, whose elements are the combinatorial objects called *vertices* and represent the fundamental units of $\Sigma$. $E$ is the set of relations, whose elements encode the relationships among the vertices, and are called the *edges*. These elements are denoted by $e$ and are the members of the power set of $V$, $\mathcal{P}(V)$, the set of all subsets of $V$ including the empty set $\emptyset$ and $V$ itself. That is, the elements of $E$ are finite subsets of $V$ of arbitrary cardinalities over $\mathbb{N}$, and $E$ thus comprises the set of relations of arbitrary arity on $V$. The set $C$ comprises the constants of the system, and may be considered as the set of all 0-ary or nullary maps on $\theta_\Sigma$. If there is no possibility of confusion and no particular need, we may not explicitly mention the set $C$ henceforth. It must be noted that the definition allows modelling all extant relations in a system, including as also going beyond the binary relationships.

For our entire discussion in what follows, we shall assume the sets $V$ and $E$ to be finite, and write $|V| = m$, $m \in \mathbb{N}$, etc, where $|V|$ denotes the number of elements in the set $V$.

Realisations of the above definition of system is instantiated by a myriad of both abstract as well as real-world examples. An example of a mathematical system is that of the natural numbers, written as $\theta_\mathbb{N} = (\mathbb{N}, <, S, +, \cdot, 0)$, where $\mathbb{N}$ is the non-empty set of elements (symbols) called numerals, $<$ is a binary relation symbol, $S$ is a unary mapping symbol, while $+$ and $\cdot$ are binary mapping symbols, often referred to as addition and multiplication, respectively, and $0$ is a constant symbol in $\mathbb{N}$. An ecosystem serves as an example of a physical system, and can be written as $\theta_E = (V_E, E_E, C_E)$, where $V_E$ is the non-empty set of all the components of a given ecosystem, $E_E$ is the set of all relations that operate among the components within the ecosystem, and $C_E$ comprises the ecological or environmental constants in the set $V_E$.

Henceforth, we shall denote a system as $\Sigma$ instead of $\theta_\Sigma$, unless stated otherwise, and when there is no possibility of confusion between the structure of the system and the system itself. Next, we define the concept of interaction.

*Interaction*

*Definition*: Let $\Sigma = (V, E, C)$ be a system. Each element of $E$ is called an *interaction* in $\Sigma$.

Thus, an interaction $e = \{v_0, \ldots, v_{k-1}\}$ containing $k$ number of vertices, $k \in \mathbb{N}$, $k \leq |V|$, is an element of the power set of $V$: $e \in \mathcal{P}(V)$. Therefore, $\forall e \in E, e \in \mathcal{P}(V)$, implying that $E \subseteq \mathcal{P}(V)$.

*Definition*: Given the above system, The *order* of an interaction $e \in E$ with $|e| = k$ is defined to be $k - 1$, and $e$ is called a $k - 1$ interaction.

Table1 lists some instances of the above definitions, with $n \in \mathbb{N}$:

**Table 1**

| Interaction $e$ | $|e|$ | Order of interaction $n$ |
|---|---|---|
| $e = \{v_0\}$ | 1 | 0 (0-interaction) |
| $e = \{v_0, v_1\}$ | 2 | 1 (1-interaction) |
| $e = \{v_0, v_1, v_2\}$ | 3 | 2 (2-interaction) |
| $\vdots$ | | |
| $e = \{v_0, \ldots, v_k\}$ | $k + 1$ | $k$ ($k$-interaction) |

Interactions are classified as *higher-order* if $|e| = k \geq 3$, and as *lower-order* if $|e| = k \leq 2$ [21].

Since each element of $E$ is an interaction in the system, therefore the *n*- ary relations that comprise $E$ encode all interactions in $\Sigma$. Owing to this fact, we shall use the term relation interchangeably with the term interaction that this relation encodes, in rest of this article. Further, if a system $\Sigma$ is such that given any interaction $e \in E$, all non-empty subsets $e'$ of $e$ also belong to $E$, then the system is said to exhibit subset dependency [22].

Having given a combinatorial description of a system comprising *n*-ary interactions, we note the identification of the term edge with the term interaction, though we shall retain both these terms throughout the work. We next proceed to represent such a system combinatorially, for which we need to introduce some objects from algebraic topology in the following section.

*Simplicial Complex*

In this section, we define and describe a combinatorial object known as simplicial complex in classical algebraic topology [23, 24]. We argue that a simplicial complex is equipped with an appropriate structure to capture interactions of all possible (finite) orders between the fundamental units of a given system which exhibits subset dependency, and constitutes an eloquent and powerful mathematical framework enabling detailed structural analysis of the resultant representation of the system. A rich body of literature has been devoted to the study of simplicial complexes, and the definitions that follow this paragraph are essentially found in the literature [25 - 46]. In this work we are interested to develop a combinatorial description of SESs and give a combinatorial definition of SEN, and hence shall remain primarily concerned with the concept of an abstract simplicial complex, while making only a passing reference to its geometric counterpart.

*Definition*: Let $V$ be a non-empty, finite set with $|V| = n + 1, n \in \mathbb{N}$, whose elements are the vertices, and are denoted by $v_i$, $i = 0, \ldots, n$. Any member of the power set of $V$, $\mathcal{P}(V)$, is a combinatorial object called a *simplex* over $V$ (or, a simplex).

For the rest of our discussion in this work, we shall assume that all simplices (also, simplexes) are non-empty. Thus, by a simplex we shall understand a non-empty, finite subset of $V$ that represents and characterises an interaction in $V$.

The *dimension* of a simplex $\sigma$ is defined as $dim(\sigma) = |\sigma| - 1$. Often, $\sigma$ called a *k*- simplex if $dim(\sigma) = k, k \in \mathbb{N}, k < n$. Any subset of a simplex $\sigma$ is called a *face* of $\sigma$. A 0-face (0-dimensional subset of a simplex) is called a vertex and written as $\{v_i\}$, and a 1-face is called an edge, written as $\{v_i, v_j\}$, with the indices belonging to some index set.

It is important to note that we have used the conventionally geometric terms such as edge and face in the above (and in what follows) only to provide names to the corresponding combinatorial objects, without supposing any geometry such as the presence of a Euclidean spacetime, for example, as a backdrop. We next define a collection of simplices as another combinatorial object, which will serve as the main building block for defining SEN.

*Definition*: Given a set $V$ as in the above definition, a family $\Delta$ of simplices is called a *simplicial complex* if it is closed under inclusion, that is, under taking of (finite, non-empty) subsets.

Though the collection $\Delta$ may either be finite or infinite, we shall assume it to be finite for our discussion. The above definition means that the following conditions hold for $\Delta$:

(1) $\forall v_i \in V, \{v_i\} \in \Delta, \ i = 0, \dots, n$. That is, for all vertices, the face $\{v_i\}$ must belong to $\Delta$.
(2) $\forall \sigma \in \Delta, (\tau \subset \sigma \Rightarrow \tau \in \Delta)$. That is, every face of a simplex in $\Delta$ must also belong to $\Delta$.

These conditions further imply that two arbitrary simplices in the simplicial complex $\Delta$ are either disjoint, or they intersect as a face of $\Delta$. That is, $\forall \sigma, \tau \in \Delta \Rightarrow (i) \ \sigma \cap \tau = \emptyset$, or $(ii) \ \sigma \cap \tau \in \Delta$. Often, the object $\Delta$ defined above is also referred to as an abstract simplicial complex.

Let $V = \{v_0, \dots, v_n\}$. Let $\mathcal{P}(V)$ be written as $\langle v_0, \dots v_n \rangle$. Thus according to the above definition, if $\sigma$ is an *n*-simplex over the vertices of $V$, then the simplicial complex $\Delta$ composed of $\sigma$ itself and the rest of its $d$ – dimensional faces ($d < n$) is written as $\Delta = \langle v_0, \dots v_n \rangle$.

For example, given a set $V = \{v_i, v_j, v_k\}$, the simplicial complex $\Delta$ comprising the 2-simplex $\{v_i, v_j, v_k\}$ and all its 0- and 1- dimensional faces is written as $\Delta = \{\{v_i\}, \{v_j\}, \{v_k\}, \{v_i, v_j\}, \{v_i, v_k\}, \{v_j, v_k\}, \{v_i, v_j, v_k\}\} = \langle v_i, v_j, v_k \rangle$. It may be noted that $\Delta$ does not include the empty simplex $\emptyset$ as per our assumption, while some authors include the empty simplex too while describing a simplicial complex [47].

The *dimension* of a simplicial complex $\Delta$, denoted by $dim(\Delta)$, is defined to be $r \geq 0$ where $r$ is the largest natural number such that $\Delta$ contains an *r*-simplex. If $dim(\Delta) = d$, then every face of dimension *d*-1 is called a *facet*.

*Definition*: Let $\Delta$ be a simplicial complex. The *boundary* of $\Delta$, denoted by $\partial \Delta$, is the set of all the faces of $\Delta$. $\partial \Delta := \{\tau | \tau \subset \sigma \in \Delta\}$, where $\sigma$ is a *k*-simplex of $\Delta$.

The above definition implies that for a given *k*-simplex in $\Delta$, every face that is a (*k*-1)- simplex is a boundary face or a facet. Further,

*Definition*: Let Δ be a simplicial complex on a set *V*. Let $p \in \mathbb{N}$. The simplicial complex given by $\Omega := \{\sigma \in \Delta | dim(\sigma) \leq p\}$ is defined to be the *p-skeleton* of Δ.

Ω is thus the collection of all faces of Δ, that have a dimension at most *p*.

For example, a simple graph $G = (V, E)$ with a finite, non-empty vertex set *V* and an edge set *E* is a 1-skeleton comprising all 0-simplices (vertices) and 1-simplices (links) of a simplicial complex Δ with $dim(\Delta) \geq 1$. Therefore, all simple graphs are simplicial complexes of dimension at most 1.

As has been remarked earlier, a simplicial complex is a combinatorial object: it is a collection of simplices satisfying conditions (1) and (2). Every simplicial complex Δ, however, corresponds to a geometric object which is a subspace of the *m*-dimensional Euclidean space $\mathbb{R}^m$ via a mapping $\varphi: V \to \mathbb{R}^m$, where *V* is the vertex set on which Δ is defined. The geometric realisation of Δ with respect to the mapping $\varphi$ is defined to be the set $|\Delta|_\varphi = \bigcup_{\sigma \in \Delta} |\sigma|_\varphi$. The object on the right hand side of the equality sign is the geometric simplicial complex, and is a subspace of $\mathbb{R}^m$, while each simplex participating in the right hand side union is the convex hull of its vertices [23, 24, 37].

**Methodology**

In order to apply the foregoing definitions of the combinatorial objects and discussions to describe a formal representation of a SES as a SEN, it is necessary to formulate the SES as a system in accordance with the definition of structure. Once this is achieved, the interactions of all orders – both lower as well as higher - can be encoded and represented as simplices and their collections as simplicial complexes of appropriate finite dimensions respectively, and the SES can be represented by a corresponding SEN that encodes the interactions in the SES. To obtain such a formulation, we propose a definition of SES from a structural perspective in the present section, followed by a definition of SEN that includes all possible orders of interactions in its structure. To argue further, we assume in this work that the focal SES exhibits subset dependency, and that the term system will henceforth mean such a system.

The above definitions of the combinatorial objects make the following observations immediate: the identification of the universe of the system with the vertex set of the simplicial complex, and the interactions of various orders in the system with simplices of various dimensions belonging to this simplicial complex. The boundary of the system is represented by the boundary of the corresponding simplicial complex, which is a collection of faces of the simplicial complex. In this formalism, all lower-order interactions in the system are encoded as 1-simplexes or 0-simplexes in the complex. However, the system may comprise higher-order interactions, each encoded as a $k$-simplex ($k \geq 2$) in the complex, which the graph theoretic pairwise interaction formalism fails to capture and represent. As a best approximation to capture the higher-order interactions, graph theoretic models of SEN, out of necessity, reduce the dimensions of the all the simplices in the complex to at most one thereby losing important structural information, and yield a simplicial complex that is actually a 1-skeleton of the higher-order SEN representing the real SES.

In order to arrive at a structural definition for SES, we premise that the fundamental ecological units of interest are the different ecological entities comprising an ecosystem and constitute the universal set $V_E$ of the constituent ecological system, while the fundamental social units are the social entities comprising the set $V_S$, the universe of the constituent social system, such that $V_S \cap V_E = \emptyset$. The last condition asserts that the social and the ecological units are separate entities, and no ecological unit is a social unit and vice-versa. We further assume that the social and ecological systems are autonomous systems and exist independently of one another, with their structures given by $\theta_S = (V_S, E_S, C_S)$ and $\theta_E = (V_E, E_E, C_E)$, respectively. The sets $E_S$ and $E_E$ comprise the social and the ecological relations respectively defined on the sets $V_S$ and $V_E$, with $E_S \cap E_E = \emptyset$. The last condition expresses the fact that the social interactions (recall, interactions encode relations) are distinct and different from the ecological interactions, and vice-versa.

*Definition:* A *social-ecological system* (SES) is an object with its structure given by $\theta_{SES} = (V, E, C)$, where the universe is the vertex set $V \neq \emptyset$ and comprises all the social-ecological units with $V = V_S \sqcup V_E$, the disjoint union of the sets $V_S \neq \emptyset$ and $V_E \neq \emptyset$. The set $E$ comprises the interactions of $\theta_{SES}$, and each interaction corresponds to an element of the power set of $V$, $\mathcal{P}(V)$. The set $C$ comprises the constants in the given SES, which are the various non-negotiables or constraints present in the system under study.

According to the foregoing definition, each element $e \in \mathcal{P}(V)$ of arbitrary cardinality, that is, an interaction of arbitrary order, comprising the collection $E$ is encoded and represented by a simplex of an appropriate dimension. Together with conditions (1) and (2), $E$ constitutes a simplicial complex of dimension $d > 0$ over $V$.

Let the environment be defined by the structure $\theta_E = (X)$ where the universe $X$ comprises all biotic and abiotic components of the environment. It may be noted that by the definition of environment, no interactions between its biotic and abiotic components are included, and hence the structure comprises only the non-empty set $X$, while the other components of the structure are empty sets, and hence are omitted. Also, since both the social as well as ecological entities are members of the environment, therefore $V_S \subseteq X, V_E \subseteq X$, with $X = V_S \sqcup V_E$. It is important to note that $E_S \sqcup E_E \subseteq E$, which captures the fact that the set $E$ of relations in the SES contains all relations of all arity and hence interactions of all orders that exist at various scales and levels of integration within the social and the ecological components individually. These interactions could be purely social-social and ecological-ecological types, but may also encode relations and interdependencies at various scales and levels of integration that exist between the social and the ecological components as cross-level relations, represented typically as social-ecological interactions [7, 9, 13].

Further, since $V = V_S \sqcup V_E \subseteq X$, there exists an inclusion map $f: V \hookrightarrow X$ defined by $v \mapsto f(v) = v \in X$. This completes the inclusion of the universe $V$ of the SES into the universe of the environment. The inclusion map embeds the focal SES within the environment. Further, we have $= V \sqcup V^c$, where $V^c = \{x | x \in X, x \notin V\}$ is the complement of the vertex set $V$ and comprises all those elements of the environment that are not included in the SES under study, and hence are external to the system, thus describing the boundary of the system.

The interacting units of the social and the ecological systems individually and independently of each other form networks of interactions between their respective fundamental units. A SEN as a formal representation of the corresponding SES comprises a collection of such networks, encoding the cross-level, inter-system interactions between the social units on the one hand and the ecological units on the other. As we have already noted, these networks do not necessarily encode only the dyadic, pairwise interactions in the system to be represented exclusively by graphs, but may encode higher-order interactions present in the systems, represented by simplicial complexes of appropriate dimensions, particularly so when group interactions with group size larger than two are being represented.

A SEN, therefore, must be understood as an empirical relational structure that captures the order-wise heterogeneous interactions within the SES, and comprises a collection of non-empty sets of various cardinalities representing the different orders of interactions that encode the $n$-ary relations that exist in the system over the set of its units, with $n \geq 1$. To formalise this concept, we define a SEN as a relational structure as follows:

*Definition:* A social-ecological network (SEN) is a combinatorial object whose structure is given by $N_{SES} = (\Delta, \Lambda)$ along with an algorithm $A$ such that for $\Lambda \neq \emptyset$, $\alpha \in \Lambda \subset \mathbb{N}$, $\Delta$ is a simplicial complex over $V$, the universe of an SES, with $\dim(\Delta) \geq 1$ given by the algorithm $A(\alpha)$. A network is called static if the temporal component $\Lambda$ is singleton, otherwise the network is a dynamic network.

The above definition implies that in the instance of an SES admitting interactions at most of order 1 (that is, only lower-order, at most binary interactions), the corresponding network $N_{SES}$ at every step indexed by $\alpha$ can be identified with a 1-skeleton of $\Delta$. Such a network is often called the underlying graph of $\Delta$, and is given by the structure $G = (V, E)$, consisting of the vertex set $V$ a non-empty set of finite abstract combinatorial objects called vertices, and a family $E$ of two-element subsets of $V$: $E \subseteq \mathcal{P}(V)$ such that $\forall e \in E, |e| = 2$, that is, each $e$ is a 1-simplex called an edge of $G$, as mentioned earlier. It is this mathematical object that gets referred to as SEN conventionally. Further, the definition of SEN implies that a given SES is uniquely represented by a SEN with respect to a given relation among the units of the SES. That is, for every $R_i$, $i = 1, \ldots, m$ in the structure of the SES, there exists a SEN that correspondingly models the SES with respect to that specific $R_i$. Due to the overarching nature of its definition, we shall refer to the defined SEN as higher-order network, while it would represent even the lower-order interactions that may be present in the system.

A rider is in order: even when a SES does not exhibit subset dependency, still the interactions in the system can be modelled using a corresponding SEN, described as a simplicial complex. In such a case, the system is represented by the facets of the complex and the information they encode, and the information about lower-order simplices is often partially lost. This cost, however is compensated by the gains from a collection of very powerful analytical methods from simplicial topology, homology and simplicial geometry that could provide very useful insight into the topology, geometry and dynamics of the system [46].

**An illustration**

We conclude the work with an illustration of the discussions above using a toy model of a SES and a corresponding SEN to specifically highlight the structural differences between the graph-theoretic SEN and the higher-order SEN defined in this work. As our example, we consider the case of a village forest in the tiger bearing Central India-Eastern Ghat landscape of India as the focal SES. This village named Saigata is located in Brahmapuri taluka of Chandrapur district in the state of Maharashtra, in Central India, in proximity to the Tadoba-Andhari Tiger reserve [48].

The villagers of Saigata, the resource users (henceforth the resource community), have shared a historical relationship with the village forest, in terms of harvesting and appropriation of forest resource. The forest was open access till the mid 1970's, and thus the forest resource was a free good, freely available to the village populace (the resource community) comprising the social system of the SES. The consequent over-consumption of the forest resource reduced the forest from thick-canopy forest to degraded land, with an accompanying loss of biodiversity and related ecosystem services. During the mid-70's, an ecological awareness coupled with cognition about their immediate future among the resource community was responsible to catalyse concerted efforts by the community to halt any further degradation of the forest, simultaneously providing the impetus to transfer the forest from open access to common property regime. This impetus led to the emergence of an informal village institution, the Van Sanrakshan Samiti (Forest Protection Committee, FPC) in 1979, with its general body, the Gramsabha, comprising the adults of each family of the village as members. The mandate of the Gramsabha was to govern the forest use through defining and evolving rules for appropriation of forest resources by the community members (round the clock community vigilance reduced the possibility of resource appropriation by any outsider to practically zero). The rules were arrived at through the integration of ideas and experiences from the local ecological knowledge of the resource community. The evolved set of rules was characterised by its being adaptive and dynamic in nature. The consequence of this participatory and collaborative process of evolution of appropriation rules was the regeneration of a resplendent village forest harbouring abundant biodiversity even occasionally visited by tigers, thereby vouching for a spectacular conservation and management success of the resource community. The key factor for this success was the community's scrupulously adhering to the adaptively evolved rules of appropriation by the FPC, without the presence of any external agent (like the state Forest Department or any other government agency) or a central controlling authority to enforce the rules.

In this toy model, we consider the adult members of the resource community to constitute the set $V_S$, while the set $E_S$ represents the social relationships among them, and these two objects comprise the social system, while the ecological entities of the village forest comprise the set $V_E$ and together with the ecological interactions within them represented in $E_E$, comprise the ecological system. The social-ecological interdependency due the historical relationship between the resource community and the village forest comprises the focal SES, with its vertex set given by $V$ and its set of interactions given as $E$ as in the definition of SES. The interdependency in the SES induces the relation of collaboration for addressing the forest resource management challenges faced by the community by evolving a set of appropriation rules for the sustainable management of the Saigata common. This relation over the set $V$ is explicitly represented as interactions of various orders in the corresponding SEN.

For modelling the SEN, let us assume that the set of rules referred to above was evolved through an algorithm requiring iterative interactions among the members of the Garamsabha. At every step of iteration, the participants were required to interact as a group starting with a group size of two at the initial step, with the output at this step being a collection of rules evolved through pairwise interactions. At each step of iteration, the sets of rules evolved by the individual groups in the previous step were the inputs, with a provision that all such sets need not be qualitatively distinct. Further, at each step, a member was required to be added to each group from another group, thus increasing the number of interacting members in each group by one per step. The groups were tasked to ideate as needed, and emerge a set of rules at each step by modifying and filtering the set emergent at the previous step. Finally when the group size reached its maximum with all members interacting as a single group, the emergent set of codes of behaviour was declared as the set of rules for appropriation of the forest resources, mandated to be adhered to by the resource community. Also, all members were required to mandatorily participate in the interactions at every iterative step, to ensure that the process is democratic, participatory and collaborative.

For the sake of simplicity, let us assume that there are a total of five different participants comprising the set $V$ in the definition of SEN, denoted as $V = \{v_i, v_j, v_k, v_l, v_m\}$, such that $v_i \neq v_j \neq v_k \neq v_l \neq v_m$. These members evolve the rules through iterative interactions among themselves within the SES, where for every step, interactions represented by the simplex $\sigma \in \mathcal{P}(V)$ encode and represent the relation of collaboration in evolving the set of appropriation rules over the set $V$. We tabulate the summary of each step of the iteration as follows, with combination of the indices over the set $\{i, j, k, l, m\}$:

**Table 2**

| Step | Input | Interaction as simplex | Order of interaction | Output | Result |
|---|---|---|---|---|---|
| 1 | 5 0-simplices | $\sigma = \{v_i, v_j\}$ $dim(\sigma) = 1$ | Lower-order | 10 1-simplices | Ten interactions with group size two each, emerging ten sets of rules. |
| 2 | 10 1-simplices | $\sigma = \{v_i, v_j, v_k\}$ $dim(\sigma) = 2$ | Higher-order | 10 2-simplices | Ten interactions with group size three each, emerging ten sets of rules. |
| 3 | 10 2-simplics | $\sigma = \{v_i, v_j, v_k, v_l\}$ $dim(\sigma) = 3$ | Higher-order | 5 3-simplices | Five interactions with group size four each, emerging five sets of rules. |
| 4 | 5 3-simplices | $\sigma = \{v_i, v_j, v_k, v_l, v_m\}$ $dim(\sigma) = 4$ | Higher-order | 1 4-simplex | One interaction with group size five, emerging one set of rules. |

As could be observed from Table 2, step 1 of the interaction results into $K_5$, the complete graph on 5 vertices. In the perspective of conventional graph-theoretic network, the interactions beyond step 1, say, at step 2, though outputs higher-order combinatorial object, would still be represented as a set of ten 1-simplices each representing a pairwise (two-person) interaction given by $\binom{5}{2}$ instead of the ten 2-simplices given by $\binom{5}{3}$, each representing a three-person interaction as in the table. Similarly, outputs of each successive step would be represented as cliques of varying sizes with each iteration resulting in a subgraph of $K_5$, encoding only pairwise interactions, with step 4 again producing a copy of $K_5$, isomorphic to and indistinguishable from the $K_5$ of step1. Though the real system relations encoded in these two steps are combinatorially entirely different, with step1 encoding pairwise interactions and step 4 encoding non-pairwise interactions, yet the respective outputs fail to capture the two very different real-world situations. This major drawback in modelling the network with exclusively pairwise interactions makes the resultant network only a coarse approximation to the real SES. While this shortcoming is owing to the structural definition of graph as was noted previously, the same is duly overcome using the simplices of appropriate dimensions to represent the higher-order interactions capturing the co-activity (that is, many-person interactions) of the participants in group sizes larger than two, as shown in Table 2, thus furnishing a model that provides a better approximation to the reality.

The structure obtained at step 4 of Table 2 is the final output of the iteration with the emergent code of behaviour as the evolved set of rules. We now use this output structure to obtain the SEN that encodes and represents this toy model of the SES. This SEN is given by $N_{SES} = (\Delta, \Lambda)$ along with an algorithm $A$ such that for $\Lambda \neq \emptyset$, $\alpha \in \{1,2,3,4\} \subset \mathbb{N}$, $\Delta$ is a simplicial complex over $V = \{v_i, v_j, v_k, v_l, v_m\}$ and is given by $\Delta = \langle v_i, v_j, v_k, v_l, v_m \rangle$ comprising all the thirty one simplices of Table 2, with $dim(\Delta) = 4$ and $\Delta \subset \mathcal{P}(V)$.

The fundamental structural differences between pairwise interactions and higher-order interactions discussed in the above illustration limits itself only to combinatorial aspects due to the stated objective of this work. However, it must be mentioned that the simplicial complex structure of the SEN induced by the higher-order interactions in the SES can lead to very interesting, rich and deep insights into the algebraic, topological and geometric structures of the SEN, which may then be studied to infer inherent structural properties and details of the complex SES that this SEN represents. Though any elaboration on this point is outside the scope of this work, the interested reader may refer to [42, 44, 46, 49, 50, 51].

**Discussion and Conclusion**

We specifically refer to the simplicial object defined as SEN in the foregoing as higher-order network, where the non-pairwise interactions are included in the representation of the corresponding SES in addition to the pairwise ones. As is clear, only those SESs that admit solely lower-order interactions at most up to order 1, that is, admit at most pairwise interactions among its fundamental units, are the ones that are represented entirely by the graph theoretic (conventional) SENs. Such instantiations are the special cases of our definition of SEN, where by setting $dim(\Delta) = 1$ we recover the conventional, graph theoretic network representing lower-order interactions in the SES.

In conclusion, this work argues that whenever a SES admits interactions beyond the pairwise ones and exhibits subset dependency, the SEN representing this system should be structurally generalised from a graph to a simplicial complex as proposed here, to include the non-pairwise, higher-order interactions and thus to obtain a more realistic representation of the SES. Such a representation is expected to be of significant modelling value, even as it would distinguish between the algebraic objects associated with the complex such as the matrices and tensors, as well as the measures and indices pertaining to different orders of interactions in the system and the corresponding scales of the network - local, mesoscale and universal [42, 44, 46]. With such an array of structural indices across different scales, the higher-order network representation would offer a deeper insight into the complex behaviour of the SES that this simplicial network represents, compared to the conventional, graph theory-based network. Representing any other SES by the structure of a graph would effectively amount to loss of information about the non-binary interactions existing in the system, and hence the structure of the SES that is being represented. Evidently, such SEN representations would result into an erroneous and incomplete insight into the structure and behaviour of the SES modelled. Based on the definition proposed in this work, the modelling framework of formally representing a SES as a SEN holds the promise of facilitating a significantly detailed insight into the structure, and hence the interactions responsible for the inherent complexity and the complex behaviour of the SES being studied.

## References

[1] Levin, S.A., 1998. Ecosystems and the Biosphere as Complex Adaptive Systems. Ecosystems 1, 431-436.

[2] Levin, S.A., Xepapadeas, T., Crepin, A-S., Norberg, J., De Zeeuw, A., Folke, C., Hughes, T., Arrow, K., Barrett, S., Daily, G., Ehrlich, P., Kautsky, N., Maler, K.G., Polasky, S., Troell, M.,Vincent, J.R., Walker, B., 2012. Social-ecological systems as complex adaptive systems: modeling and policy implications. Environment and Development Economics 18, 111–132.

[3] Levin, S. A., Clark, W.C. (Eds.), 2010. Toward a Science of Sustainability. Center for International Development Working Papers 196, John F. Kennedy School of Government, Harvard University.

[4] Bodin, O., Crona, B.I., 2009. The role of social networks in natural resource governance: What relational patterns make a difference? Global Environmental Change 19, 366–374.

[5] Bodin, O., Prell, C. (Eds.), 2011. Social Networks and Natural Resource Management: Uncovering the Social Fabric of Environmental Governance. Cambridge University Press.

[6] Gonzale`s, R., Parrott, L., 2012. Network theory in the assessment of the sustainability of social–ecological systems. Geography Compass 6/2, 76–88.

[7] Bodin, O., Tengo, M., 2012. Disentangling intangible socialecological systems Glob. Environ. Change. 22, 430–9.

[8] Bodin, O., Crona, B., Thyresson, M., Golz, A-L., Tengö,M., 2014. Conservation success as a function of good alignment of social and ecological structures and processes. Conserv. Biol. 28, 1371–9.

[9] Bodin, Ö., Robins, G., McAllister, R. R. J. , Guerrero, A. , Crona, B., Tengö, M., Lubell, M., 2016. Theorizing benefits and constraints in collaborative environmental governance: a transdisciplinary social-ecological network approach for empirical investigations. Ecology and Society 21(1), 40.

[10] Baggio, J.A., BurnSilver, S.B., Arenas, A., Magdanz, J.S., Kofinas, G.P., De Domenico, M., 2016. Multiplex social ecological network analysis reveals how social changes affect community robustness more than resource depletion Proc. Nat. Acad. Sci. U.S.A. 113 (48), 13708–13713.

[11] Bodin, O., Barnes, M.L., McAllister, R.R.J., Rocha, J.C., Guerrero, A.M., 2017. Social-Ecological network approaches in interdisciplinary research: a response to Bohan et al and Dee et al. Trends Ecol. Evol. 32, 547–9.

[12] Barnes, M. L., Bodin, Ö. ., Guerrero, A. M., McAllister, R. J., Alexander, S. M., Robins, G., 2017. The social structural foundations of adaptation and transformation in social–ecological systems. Ecology and Society 22(4), 16.

[13] Bodin, O., Alexander, S.M., Baggio, J., Barnes, M.L., Berardo, R., Cumming, G.S., Dee, L.E., Fischer, A.P., Fischer, M., Mancilla Garcia, M., Guerrro. A.M., Hileman, J., Ingold, K., Matous, P., Morrison, T.H., Nohrstedt, D., Pittman, J., Robins, G., Sayles, J.S., 2019. Improving network approaches to the study of complex social–ecological interdependencies Nature Sustainability. 2, 551–559.

[14] Sayles, J. S., Mancilla Garcia, M., Hamilton, M., Alexander, S.M., Baggio, J.A., Fischer, A.P., Ingold, K., Meredith, G.R., Pittman, J., 2019. Social-ecological network analysis for sustainability sciences: a systematic review and innovative research agenda for the future. Environ. Res. Lett. 14, 093003.

[15] Hamilton, M., Paige Fischer, A., Ager, A., 2019. A social-ecological network approach for understanding wildfire risk governance. Global Environmental Change 54, 113 – 123.

[16] Kluger, L.C., Gorris, P., Kochalski, S., Mueller, M.S., Romagnoni, G., 2020. Studying human–nature relationships through a network lens: a systematic review. People and Nature. 2, 1100–1116.

[17] Metzger, J. P., Fidelman, P.,Sattler, C., Schröter, B., Maron, M., Eigenbrod, F., Fortin, M-J., Hohlenwerger, C., Rhodes, J. R., 2021. Connecting governance interventions to ecosystem services provision: a social-ecological network approach. People and Nature 3, 266–280.

[18] van Dalen, D., 2008. Logic and Structure. Springer.

[19] Shoenfield, J.R., 1967. Mathematical Logic. Addison-Wesley.

[20] Srivastava, S. M., 2008. A Course on Mathematical Logic. Springer.

[21] Battiston, F., Cencetti, G., Iacopini, I., Latora, V., Lucas, M., Patania, A., Young, J-G., Petri G., 2020. Networks beyond pairwise interactions: structure and dynamics. Phys. Rep. 874, 1–92.

[22] Torres, L., Blevins, A.S., Bassett, D., Eliassi-Rad, T., 2021. The why, how, and when of representations for complex systems. SIAM Rev. 63 (3),435–485.

[23] Hatcher, A., 2002. Algebraic Topology. Cambridge University Press.

[24] Deo, S., 2018. Algebraic Topology: A Primer. Springer.

[25] Muhammad, A., Egerstedt, M., 2006. Control using higher order Laplacians in network topologies. In: Proceedings of 17th International Symposium on Mathematical Theory of Networks and Systems, pp. 1024–1038.

[26] Muhammad, A ., Jadbabaie, A ., 2007. Decentralized computation of homology groups in networks by gossip. In: American Control Conference, 2007. ACC'07. IEEE, pp. 3438–3443.

[27] Maletic, S., Rajkovic, M., 2009. Simplicial complex of opinions on scale-free networks. In: Complex Networks. Springer, Berlin, Heidelberg, pp. 127–134.

[28] Horak, D., Maletic, S., Rajkovic, M., 2009. Persistent homology of complex networks. J. Stat. Mech. 2009 (03), P03034 .

[29] Edelsbrunner, H., Harer, J.L., 2010. Computational Topology: An Introduction, American Mathematical Society.

[30] Maletic, S., Rajkovic, M., 2012. Combinatorial Laplacian and entropy of simplicial complexes associated with complex networks. Eur. Phys. J. Spec. Top. 212 (1), 77–97.

[31] Petri, G., Scolamiero, M., Donato, I., Vaccarino, F., 2013. Topological strata of weighted complex networks. PLoS One 8 (6), e66506.

[32] Maletic, S., Rajkovic, M., 2014. Consensus formation on a simplicial complex of opinions. Physica A 397, 111–120.

[33] Edelsbrunner, H., 2014. A Short Course in Computational Geometry and Topology. Springer.

[34] Ghrist, R. W., 2014. Elementary Applied Topology. Createspace.

[35] Petri, G., Expert, P., Turkheimer, F., Carhart-Harris, R., Nutt, D., Hellyer, P.J., Vaccarino, F., 2014. Homological scaffolds of brain functional networks. Jour. Roy. Soc. Interface 11 (101), 20140873.

[36] Kahle, M., 2014. Topology of random simplicial complexes: a survey. AMS Contemporary Mathematics, 260, 201 – 222.

[37] Nanda, V., Sazdanovic, R., 2014. Simplicial models and topological inference in biological systems. Discrete and Topological Models in Molecular Biology, Natural Computing Series pp. 109 -141, Springer.

[38] Bianconi, G., 2015. Interdisciplinary and physics challenges of network theory. Europhys. Lett. 111 (5), 56001.

[39] Bianconi, G., Rahmede, C., 2016. Network geometry with flavor: from complexity to quantum geometry. Phys. Rev. E 93 (3), 032315.

[40] Sizemore, A., Giusti, C., Bassett, D.S., 2016. Classification of weighted networks through mesoscale homological features. J. Complex Netw. cnw013 .

[41] Giusti, C., Ghrist, D., Bassett, D. S., 2016. Two's a company, three (or more) is a simplex. Jour. Comp. Neurosciences 41(1), 1 – 14.

[42] Courtney, O.T., Bianconi, G., 2016. Generalized network structures: The configuration model and the canonical ensemble of simplicial complexes. Phys. Rev. E 93, 062311.

[43] Courtney, O.T., Bianconi, G., 2017. Weighted growing simplicial complexes. Phys. Rev. E 95, 062301.

[44] Estrada, E., Ross, G.J., 2018. Centralities in simplicial complexes. Applications to protein interaction networks. Jour. Theor. Biol. 438, 46 – 60.

[45] Icaponi, I., Petri, G., Barrat, A., Latora, V., 2019. Simplicial models of social contagion. Nature Communications 10 (10), 1 – 9.

[46] Bianconi, G., 2021. Higher Order Networks: An Introduction to Simplicial Complexes. Cambridge University Press.

[47] Jonsson, J., 2008. Simplicial Complexes of Graphs. Springer-Verlag.

[48] Bhattacharya, S., 2004. The evolution of cooperation at Saigata commons: A game theoretic interpretation. Working Paper No. 6-04, Ostrom Workshop, Indiana University, Bloomington.

[49] Horak, D., Jost, J., 2013. Spectra of combinatorial Laplace operators on simplicial complexes. Advances in Mathematics 244, 303–336.

[50] Torres, J. J., Bianconi, G., 2020. Simplicial complexes: higher-order spectral dimension and dynamics. J.Phys.Complex. 1, 015002.

[51] Schaub, M.T., Benson, A.R., Horn, P., Lippner, G., Jadbabaie, A., 2020. Random walks on simplicial complexes and the normalized Hodge 1-Laplacian. SIAM Rev. 62 (2), 353–391.